\documentclass[11pt,reqno]{amsart}
\usepackage{hyperref}
\usepackage{amsmath,amsfonts,amssymb}
\usepackage[latin1]{inputenc}
\usepackage{verbatim}
\usepackage{color}

\newtheorem{theorem}{Theorem}[section]
\newtheorem{lemma}[theorem]{Lemma}
\newtheorem{prop}[theorem]{Proposition}

\newtheorem{rmk}[theorem]{Remark}
\newtheorem{defn}[theorem]{Definition}

\newcounter{defn}

\newcommand{\sect}{\vspace{3mm} \setcounter{equation}{0} \setcounter{defn}{0} \section}

\newcommand{\w}[1]{\langle {#1} \rangle}

\newcommand{\pf}{\noindent {\bf Proof. \hspace{2mm}}}
\newcommand{\ef}{ \hfill $ \Box $ \vskip 3mm}
\newcommand{\be}{\begin{equation}}
\newcommand{\ee}{\end{equation}}
\newcommand{\bea}{\begin{eqnarray}}
\newcommand{\eea}{\end{eqnarray}}
\newcommand{\beas}{\begin{eqnarray*}}
\newcommand{\eeas}{\end{eqnarray*}}
\newcommand{\ep}{{\epsilon}}
\newcommand{\f}{\frac}
\newcommand{\na}{\nabla}

\newcommand{\nn}{\nonumber}

\newcommand{\bC}{{\mathbb C}}
\newcommand{\bR}{{\mathbb R}}
\newcommand{\bN}{{\mathbb N}}
\newcommand{\bZ}{{\mathbb Z}}

\newcommand{\vF}{{\mathcal F}}
\newcommand{\vG}{{\mathcal G}}

\newcommand{\vH}{{\mathcal H}}

\newcommand{\vL}{{\mathcal L}}
\newcommand{\vM}{{\mathcal M}}

\newcommand{\vT}{{\mathcal T}}

\newcommand{\vS}{{\mathcal S}}

\newcommand{\gM}{{\mathfrak M}}

\def\f{\frac}
\def\Dl{\Delta}
\def\na{\nabla}
\def\la{\lambda}

\def\Dl{\Delta}

\def\i{\infty}

\def\lt{\left}
\def\rt{\right}

\begin{document}

\title[The Kramers-Fokker-Planck equation]{
The Kramers-Fokker-Planck equation with a decaying potential in $\bR^n$,  $n \ge 4$}

\author{Xinghong Pan}
\address[X. Pan]{School of Mathematics and Key Laboratory of MIIT, Nanjing University of Aeronautics and Astronautics, Nanjing 211106, China}
\email{xinghong\_87@nuaa.edu.cn}

\author{Xue Ping WANG}
\address[X.P. Wang]{School of Mathematics and Key Laboratory of MIIT, Nanjing University of Aeronautics and Astronautics, Nanjing 211106, China}
\email{xpwang@nuaa.edu.cn}

\author{Lu ZHU}
\address[L. Zhu]{School of Mathematics, Hohai University, Nanjing 211100, China}
\email{zhulu@hhu.edu.cn}

\subjclass[2000]{35J10, 35P15, 47A55, 82C40}
\keywords{Kramers-Fokker-Planck equation, nonequilibrium statistical mechanics, large-time behavior of solutions, non-selfadjoint operators, threshold spectral analysis, quantum scattering}

\begin{abstract}

We use methods from microlocal analysis and quantum scattering to study spectral properties near the threshold zero of the Kramers-Fokker-Planck operator with a decaying potential in $\bR^n$, $n \ge 4$, and deduce the large-time behavior of solutions to the kinetic Kramers-Fokker-Planck equation.  For short-range potentials, we establish an optimal time-decay estimate in weighted $L^2$-spaces when $ n\ge 5$ is odd.  For potentials decaying like $O(|x|^{-\rho})$ for some $\rho > n-1$, we obtain, for all dimensions $n \ge 4$, a large-time expansion of the solution with the leading term given by the Maxwell-Boltzmann distribution multiplied by the factor $(4\pi t)^{-\f n 2}$ corresponding to the decay for the heat equation. These results complete those obtained in \cite{nw,w2} for dimensions $n=1$ and $3$. The same questions for $n=2$ are still open.
\end{abstract}

\maketitle

\sect{Introduction}

In this work, we study the large-time behavior of solutions to  the Kramers-Fokker-Planck (KFP, for short) equation,
also called  the Kramers euqation by physicists (\cite{kr,risc}), given by
\be\label{equation}
\f{\partial u(t)}{\partial t} + P u(t) =0,
\ee
where $u(t) = u(t; x, v)$ is the distribution function  of particles in the phase space $\bR^{2n}$ verifying some initial condition  $u(0) =u_0$, $P$ is the KFP operator
\be\label{operator}
 P=-\Delta_v+\f{1}{4}|v|^2-\f{n}{2} + v\cdot\na_x-\na V(x)\cdot\na_v, \quad (x,v) \in \bR^{2n}.
  \ee
  $V(x)$ is supposed to be  a  real-valued $C^1$ potential satisfying, for some constants $\rho$ and $ C>0$ the estimate
\be \label{ass1}
 |V(x)| + \w{x}|\nabla V(x)| \le C\w{x}^{-\rho},  \quad x\in \bR^n,
 \ee
 with $\w{x} = (1 + |x|^2)^{1/2}$.

 The large-time asymptotics of the KFP equation is initially motivated by a mathematical proof of the retour to equilibrium, a physical phenomenon predicted by the second law of thermodynamics for isolated systems. The condition {\it ``isolated systems''} is modeled in mathematics either by confining potentials $V(x) \to + \infty$ (\cite{dv, hrn, hln, lz, v}) or by bounded domains with appropriate boundary conditions (\cite{leb,n}). The typical result is that for symbol-like confining potentials with various additional conditions, the distribution function converges to the global equilibrium, which is the Maxwell-Boltzmann distribution in these cases, with an explicit convergence rate as $t\to +\infty$.

 We are interested in the KFP equation with a decaying potential ($\rho >0$) which can be used to model  non-equilibrium thermodynamic systems where particles interact with transient, localized force. A typical example is the free gas expansion with white noise.  For decaying potentials, the diffusion dominates over the drift and the irreversibility of process predicts that the solution tends to zero at large times.  A natural mathematical question is to  justify rigorously this physical prediction and to calculate the large-time asymptotics of the solution as precisely as possible. The first result we are aware of in this direction is \cite{w2}, where the second author of the present work, combining methods from microlocal analysis,  spectral analysis for non-selfadjoint operators and quantum scattering theory, obtained optimal time-decay and large-time expansion of solutions in weighted $L^2$-spaces for short-range potentials for $n=3$. Similar results are obtained in \cite{nw} for one-dimensional case for very quickly decaying potentials.  An essential difference with the KFP equation with confining potentials is that for decaying potentials, there is no global equilibrium. The Maxwell-Boltzmann distribution $\gM$ given by
\be \label{MB}
\gM(x,v)= \f{1}{(2\pi)^{\f n 4}}  e^{- \f 1 2( \f{|v|^2} 2 + V(x))}
\ee
satisfies the stationary equation $P \gM =0$, but it is only a local equilibrium if $V(x)$ is bounded. Spectral properties are also quite different in the both cases: while the spectrum of the KFP operator are purely discrete when the potential is confining and $|\nabla V(x)| \to \infty$ (\cite{hln,wxl}), those for decaying potentials contain the interval $[0,+\infty[$ as the essential spectrum (\cite{w2}). The latter suggests that the method from quantum scattering may be pertinent for decaying potentials.

Let us mention that in \cite{bdmms} the authors attempt to set up a frame-work for hypocoercivity without confinement. They only study some translation-invariant kinetic equations. When specified to the KFP equation, this means that the condition $V(x) =0$ is required.

In this work, we study the KFP equation with a decaying potentials of dimensions $n\ge 4$.
Denote $P_0$ the free KFP operator (with $V(x) =0$):
\[
P_0 = -\Delta_v+\f{1}{4}|v|^2-\f{n}{2} + v\cdot\na_x, \quad (x,v) \in \bR^{2n}.
\]
$P_0$ defined on  $\mathcal{C}_0^{\infty}(\mathbb{R}^{2n})\;$ is essentially maximally
accretive  (see Proposition 5.5, \cite{hln}). We still denote by $P_0$
its closed extension with maximal domain $D(P_0)=\{f \in L^2(\bR^{2n}); P_0 f \in L^2(\bR^{2n})\}$.  One has:
$\sigma(P_0) =[0, \infty[$. The following subellitptic estimate holds true for $P_0$ (\cite{hln}):
\be \label{subelliptic}
\|\Delta_v f\| +\|v^2 f\| + \|\w{D_x}^{\f 2 3}f\| \le C (\|P_0 f\| + \|f\|) \quad f \in D(P_0).
\ee
For $r\ge 0,s \in \bR$, we introduce accordingly the weighted Sobolev spaces:
\[
\vH^{r,s} = \{u \in \vS'(\bR^{2n}); (1 + \w{D_v}^2 + |v|^2 + \w{D_x}^{\f 2 3})^{\f r 2 } \w{x}^s u \in L^2\}.
\]
For $r <0$ and $s \in \bR$, $\vH^{r,s}$ is taken to be the dual space of $\vH^{-r,-s}$.  Let $\vL(r,s; r',s')$ be the space of
 bounded linear operators from $\vH^{r,s}$ to $\vH^{r',s'}$.
To simplify notation we denote  $\vH^s =\vH^{0,s}=L^2(\bR^{2n}, \w{x}^{2s}dxdv)$, $\vH = \vH^{0,0}= L^2(\bR^{2n})$ and $\vL(s,s') = \vL(0,s; 0,s')$.
\\

Write the full KFP operator $P$ as
\[
P= P_0 + W \quad \mbox{ with }W = -\na_x V(x) \cdot \na_v.
\]
Under the condition (\ref{ass1}) with $\rho>-1$, $\nabla_xV(x)$
tends to $0$ as $|x| \to \infty$.  (\ref{subelliptic}) implies that $W$
is relatively compact with respect to $P_0$. Therefore $P$
defined on $D(P) =D(P_0)$ is closed and maximally accretive.
Let $\sigma (P)$ (resp., $\sigma_d(P)$, $\sigma_p(P)$) denote the spectrum (resp.,
the discrete spectrum, the point spectrum) of $P$.
$\sigma_d(P)$ is the set of isolate eigenvalues with finite (algebraic) multiplicity. Let
$\sigma_{\rm ess}(P) =\sigma(P) \setminus \sigma_d(P)$ be the essential spectrum of $P$.
When $\rho >-1$, zero may be an embedded eigenvalue if $V(x)$  increases sublinearly. Under some additional conditions,
the return to equilibrium is proven in \cite{lz} with a subexpoential convergence rate.
\\

 The general theory on relatively compact perturbation (\cite{d}) for closed non-self-adjoint operators affirms that
\be
\sigma_{\rm ess} (P) = \sigma (P_0)=[0, \infty[.
\ee
and $\sigma(P) \setminus [0, \infty[ =\sigma_d(P)$ and that the discrete spectrum $\sigma_d(P)$ may accumulate towards any
point in $[0, \infty[$. The set of  thresholds, $\vT= \bN$,  of the KFP operators is composed of  eigenvalues of the harmonic oscillator
 $ -\Delta_v+\f{1}{4}|v|^2-\f{n}{2}$ which is the real part of $P_0$.  The term threshold is appropriate here, because each eigenvalue $k$ contributes the segment $[k, +\infty[$ to $\sigma(P_0)$ (see \cite{w2}). More precise spectral properties and exponential decay of eignefunctions are studied in \cite{w4}.
\\

The main results of this work are the following

\begin{theorem}\label{th5.4} Let $S(t) =e^{-tP}$, $t\ge 0$,  be the $C_0$-semigroup of contractions generated by $-P$.

 (a).   Let $n\ge 5$ be odd. Assume condition (\ref{ass1}) with $\rho>1 $. Then one has, for any $s> \f {n}2$,
\be \label{eq5.8}
\|S(t)\|_{\vL(s; -s)} \le C_s  \w{t}^{- \f n 2}, \quad t  >0.
\ee

(b). Let $n\ge 4$. Assume condition (\ref{ass1}) with $\rho >n-1$. Then  for any $s> \f n 2$, there exists some
constant $\ep>0$ such that
\be \label{eq5.7} S(t) =
\f{ 1 }{(4 \pi t)^{ \f n
2}}\w{\cdot, \gM}\gM    + O(t^{- \f n 2 -\ep})
\ee
in $\vL(s; -s)$ as $t \to + \infty$, where $\gM$ is the Maxwell-Boltzmann distribution.
\end{theorem}

\begin{rmk}
Theorem \ref{th5.4} (b) implies that if $u_0 \in \vH^{s}$, $ s > \f n 2$, and $\w{u_0, \gM} =0$, then the solution $u(t)= S(t)u_0$ of the Cauchy problem \eqref{equation} satisfies, for some positive constants $\ep$ and $C$,
\[
\|u(t)\|_{-s} \le C \w{t}^{- \f n 2 -\ep} \|u_0\|_{s}, \quad t \ge 0.
\]
By studying the smoothness of the boundary values of the resolvent near $0$, this upper-bound may be improved if the potential $V(x)$  decays faster.
\end{rmk}

The proof of Theorem \ref{th5.4} is based on the low-energy expansion of the resolvent $R(z) = (P-z)^{-1}$ near the first threshold zero. To
do this, we use the techniques of threshold spectral analysis known for Schr\"odinger operators, with a notable difference that for Schr\"odinger operators, low-energy resolvent expansions are established according to whether zero is an eigenvalue and/or a resonance (see \cite{jn,sw} for an overview and more references), for the KFP operator we do not use any such implicit spectral condition. See \cite{w2} for the notion of threshold resonances for the KFP operator. It is well-known that threshold spectral analysis is model- and dimension-dependent \cite{jn,sw}. This explains why we distinguish the cases $n\ge 5$ odd and $n \ge 4$ even, in addition to some ideas and formulas from \cite{w2}.

The remaining part of this work is organized as follows. In Section 2, we recall from \cite{w2} some formulas needed in this article. Section 3 is devoted to the expansion of the free resolvent near threshold zero in weighted Sobolev spaces. We distinguish the cases $n\ge 5$ odd and $n\ge 4$ even, because of the low-energy behaviors of Green function for the Laplacian.  In Section 4, we study  spectral properties of $P$ at threshold zero and establish the resolvent expansions. A particular attention is paid to the calculation of the first  term with non integral power in $z$ which is needed to compute explicitly  the large-time asymptotics of the semigroup $S(t)$. Theorem \ref{th5.4} is proven by representing $S(t)$ as Cauchy integral of the resolvent. As an application of Theorem \ref{th5.4}, we give in Section 6 some dispersive estimates for $S(t)$  in $L^p-L^q$ setting.

\sect{Preliminaries}

We first recall from \cite{w2} some results needed in this work. Let
\[
P_0 = -\Delta_v+\f{1}{4}|v|^2-\f{n}{2} + v\cdot\na_x, \quad (x,y) \in \bR^{2n},
\]
and $R_0(z) = (P_0-z)^{-1}$.
By the partial Fourier transform $\vF_{x\to \xi} $ in $x$-variables, the free KFP
$P_0$ is a direct integral of the family $\{\widehat{P}_0(\xi) ; \xi \in \bR^n\}$, where
\be
\widehat{P}_0(\xi) = -\Dl_v+\f{v^2}{4} -\f{n}{2}+ i v\cdot\xi.
\ee
The operator $\widehat{P}_0(\xi)$ is given by
\[
\widehat{P}_0(\xi) = \vF_{x\to \xi} P_0 \vF_{\xi \to x}^{-1} = -\Dl_v+\f{1}{4}\sum^n_{j=1}(v_j+2i\xi_j)^2-\f{n}{2}+|\xi|^2.
\]
 $\{\widehat{P}_0(\xi), \xi\in \bR^n\}$ is  a holomorphic family of type $(A)$ in sense of Kato with constant domain
 $D= D(-\Dl_v+\f{v^2}{4})$ in $L^2(\bR^n_v)$.
  Let $F_j(s)=(-1)^je^{\f{s^2}{2}}\f{d^j}{ds^j}e^{-\f{s^2}{2}}, j \in \bN,$ be  Hermite polynomials and
$$
\varphi_j(s)=(j!\sqrt{2\pi})^{-\f{1}{2}}e^{-\f{s^2}{4}}F_j(s)
$$
the normalized Hermite functions.  For $\xi \in \bR^n$ and $\alpha=(\alpha_1, \alpha_2, \cdots, \alpha_n) \in \bN^n$, define
\be
\psi_\alpha(v) = \prod_{j=1}^n\varphi_{\alpha_j}(v_j) \mbox{ and } \psi_\alpha^\xi(v) = \psi_\alpha(v + 2i \xi).
\ee
One can check (\cite{w2}) that each eigenvalue $E_\ell(\xi) = \ell +|\xi|^2$, $\ell \in \bN$,  of $\widehat{P}_0(\xi)$  is semi-simple  and
the associated Riesz projection is given by
 \be\label{projection}
 \Pi^\xi_\ell\phi=\sum_{\alpha\in \bN, |\alpha|=\ell} \langle\psi^{-\xi}_\alpha , \phi \rangle\psi^\xi_\alpha, \quad \phi \in L_v^2.
  \ee

For a temperate symbol $a(x,v; \xi, \eta)$ (\cite{hor}), denote by $a^w(x,v, D_x, D_v)$ the associated Weyl pseudo-differential operator defined by
\bea
\lefteqn{a^w(x,v,D_x,D_v)u(x,v)} \\[2mm]
& =& \f{1}{(2\pi)^{2n}} \int\int e^{i(x-x')\cdot\xi + i(v-v')\cdot\eta} a(\f{x+x'}2, \f{v+v'}2, \xi,\eta) u(x', v') dx'dv'd\xi d\eta \nonumber
\eea
for $u \in \vS(\bR^{2n}_{x,v})$.

\begin{prop}[\cite{w2}]\label{prop2.1}  Let $\ell\in \bN$ and $ \ell < a <\ell +1$ be fixed. Take  $\chi \ge 0$ and
$\chi \in C_0^\infty(\bR^n_\xi)$ with supp $\chi \subset \{\xi, |\xi|^2 \le a + 4\}$, $\chi(\xi) =1$ when $|\xi|^2 \le a +3$ and $ 0 \le \chi(\xi) \le 1$.  Then one has
\be \label{RR0}
 R_0(z) = \sum_{k=0}^\ell b_k^w(v, D_x, D_v)(-\Delta_x +k -z)^{-1} + r_\ell(z),
\ee
for  $z \in \bC$ with $\Re z <a$ and $\Im z \neq 0$. Here  $  b_k^w(v, D_x, D_v)$ is a  Weyl pseudodifferential operator with symbol $b_k$ independent of $x$ given by
 \be
 b_k(v, \xi,\eta) = \int_{\bR^n} e^{-i v'\cdot\eta/2}\left(\sum_{|\alpha|=k} \chi(\xi) \psi_\alpha( v+ v' + 2i\xi)\psi_\alpha( v- v' + 2i\xi) \right) dv'.
 \ee
 In particular,
\be
b_0(v, \xi,\eta) = 2^{\f n2}\chi(\xi)e^{-v^2 -\eta^2 + 2iv\cdot\xi + 2\xi^2}.
\ee
 $ r_\ell(z)$ is  bounded on $L^2$ and  holomorphic in $z$ with $\Re z <a$ satisfying
 the estimate
\be \label{reste}
\sup_{\Re z < a} \|r_\ell(z)\|_{\vL(\vH)} <\infty.
\ee
\end{prop}

 $b_k^w(v, D_x, D_v)$  can be written in terms of the Riesz projection $\Pi_k^{\xi}$:
\be
b_k^w(v, D_x, D_v) =\chi(D_x) \Pi_k^{D_x}.
\ee
It follows that
\be
P_0 \; b_k^w(v, D_x, D_v) = (-\Delta_x+k) b_k^w(v, D_x, D_v).
\ee

The continuity of the free resolvent in the weighted Sobolev spaces $\vH^{r,s}$ can be summarized as follows.

 \begin{lemma}[\cite{w4}] \label{lem2.3} Let $z\in \bC\setminus \bR_+$. For any $m\in \bN$,
 $ s \in \bR$ and $r\in[0, 2m]$, one has $R_0(z)^m \in \vL(-r, s; 2m-r,s)$.
 \end{lemma}

 From Proposition \ref{prop2.1} and the limiting absorption principles for Schr\"odinger operators (\cite{ag1}), one derives the following

\begin{prop}\label{prop2.3} Let $n \ge 1$. For any $s>\f 1 2$,
the boundary values of the resolvent
\[
R_0(\lambda \pm i 0) =  \lim_{ \pm \Im z >0, z \to \lambda}R_0(z)
\]
 exist
in $\vL(0,s; 2,-s)$ for $\lambda \in \bR_+\setminus \vT$ and is continuous in $\lambda$. More generally,
let $k \in \bN$ and $s> k + \f 1 2$ and $ r, r' \ge 0 $ with $r+r' \le 2(k+1)$, $R_0(\lambda
\pm i0)$ is of class $C^k$ in $\vL(-r,s; r',-s)$ for $\lambda \not\in \vT$.
\end{prop}

The following result is useful to compute explicitly the physically relevant term in large-time expansion of the solution.

\begin{prop}\label{key-prop} Let $n \ge 1$ and  $u \in \vH$.
Then
\be \label{e4.100}
\lim_{\lambda \to 0_-} \lambda R_0(\lambda)u =0, \quad \mbox{ in } \vH.
\ee
\end{prop}
\pf (\ref{RR0}) gives
\[
 R_0(\lambda) u = b_0^w(v,D_x, D_v) (-\Delta_x -\lambda)^{-1}u + r_0(\lambda)u.
\]
Clearly, $\lambda r_0(\lambda)u \to 0 $  as $ \lambda \to 0_-$, because $r_0(z)$ is uniformly bounded for $z$ near $0$ and $u\in \vH$.
For $g \in L^2_x$, Lebesgue's Dominated Convergence Theorem shows that
\[
\lim_{\lambda \to 0_-} \lambda (|\xi|^2-\lambda)^{-1}\hat{g}(\xi) =0
\]
 in $L^2_\xi$. Consequently, for $u \in \vH$,
\[
\lim_{\lambda \to 0_-} \lambda b_0^w(v,D_x, D_v) (-\Delta_x -\lambda)^{-1}u =0, \quad \mbox{ in } \vH.
\]
(\ref{e4.100}) is proven. \ef

\sect{Low-energy expansions for the free resolvent}

We consider the free resolvent $ R_0(z) $ near the first threshold $0$  when $ n\ge 4$. Its low-energy behavior is dependent on space dimension. We distinguish the cases $n\ge 5$ odd and $ n \ge 4$ even.

\subsection{$n\ge 5$ odd.} Recall that for $n\ge 3$ odd, the integral kernel for $(-\Delta_x -z)^{-1}$, $ z\not\in [0, \infty[$,  is given by
\be\label{Green}A(x, y; z) = \f i 4 \left( \frac{z^{\f 12}}{2 \pi
|x-y|}\right)^{\f n 2 -1}H_{\f n 2 -1}^{(1)}( z^{\f 12} |x-y|),
 \ee
 where $H_\nu^{(1)}$ are the modified Hankel functions (\cite{ol})
 and the branch of $z^{\f 12}$ is taken such that $\Im z^{\f 12}>0$ for $z
 \not\in [0,\infty[$. Using the properties of the Hankel functions, this integral kernel can be
expanded asymptotically in $z^{\f 12} |x-y|$ as
\begin{eqnarray} \label{kernelnodd}
A(x, y; z)  & = & \f 1{|x-y|^{n-2}}  \left\{\sum_{k=0}^{\f{n-3} 2} a_{n,2k} z^k \zeta^{2k}
  +  \sum_{k= n-2}^\infty a_{n, k} z^{\f k 2} \zeta^k \right\}
\end{eqnarray}
where $ \zeta =  |x-y|$ and $a_{n,k}$ are numerical constants
which can be computed explicitly (\cite{jn}). In particular,
$a_{n,k} =0$  for $k$ odd and $0 \le k\le \f{n-3} 2$. The above
expansion means that for $N \in \bN$ and $0 < \epsilon \le 1$, one
has
\[
A(x, y; z)   =  \f 1{|x-y|^{n-2}}  \left\{
    \sum_{k= 0}^N a_{n, k} z^{\f k 2} \zeta^k  + O((|z|^{ \f 1 2 } |\zeta|)^{ N  + \epsilon})
  \right\}.
\]

 Let $A_k$ denote the operator with integral kernel \be A_k(x, y)
=  a_{n, k} |x-y|^{k- (n-2)}, \quad k \in \bN. \ee For $0
\le k\le \f{n-3} 2$, $A_{2k}$ can be interpreted as
\[
\w{A_{2k} u, v} = \lim_{r <0, r\to 0_-}\f{1}{k!}\f{d^k}{dr^k}\w{ (-\Delta_x-r)^{-1}u, v}, \quad u, v \in \vS(\bR^n_x).
\]
  $A_{n-2}$ is the operator of rank one accociated with the constant integral kernel $a_{n,n-2}$.  By Theorem 2.2  of \cite{w0} (with $q(\theta)= 0$) on  low-energy expansions for the Laplace-Beltrami operator on conical manifolds, one has
\[
\w{x}^{-s} A_{2k}\w{x}^{-s} \in \vL(L^2(\bR^n_x)),  \qquad s>2k+1,
\quad  \forall k\in \bN.
\]
For the first few terms, one has a better result using the properties
of integral operators with homogeneous kernels (Lemma 2.1, \cite{nw}):  for $ 0 < \alpha <n$ and  for real
numbers $a, b$ with $a < \f n 2 $, $ b < \f n 2$ and $a +b =
\alpha$, the integral kernel
\[
\f{1}{|x|^{a} |x-y|^{n - \alpha} |y|^b},
\]
defines a bounded operator in $L^2(\bR^n_x)$. Applying this result
to $\alpha = 2 + 2k$, $ k =0, 1, \cdots, \f{n-3}{2}$, one deduces
that
 \be \label{A2k} \w{x}^{-s'} A_{2k}\w{x}^{-s} \in
\vL(L^2(\bR^n_x)), \qquad  k =0, 1, \cdots, \f{n-3}{2},
 \ee
for $s, s' > \f{4-n} 2 +2k$ and $ s+ s' \ge 2k +2$. Clearly,
$\w{x}^{-s'} A_{n-2}\w{x}^{-s} \in \vL(L^2(\bR^n_x))$ for $s, s' >
\f n 2$. Making use of (\ref{RR0}) with $\ell=0$, we deduce the
following expansions for the resolvent $R_0(z)$ of the free KFP
operator $P_0$.

\begin{lemma} \label{lem2.5} For $\delta>0$ small, denote $\Omega_\delta =\{z\in \bC; |z| <\delta, z\not \in \bR_+\}$. For $n \ge 5$ odd,  the following expansions hold true.

\begin{enumerate}
\item For any $s>\f{4-n} 2$, $s'>\f{4-n} 2$ and $s+s'> 2$, there
exists some $\epsilon
>0$ such that \be \label{eq2.16} R_0(z) = G_0 + O(|z|^\ep), \quad \mbox{ for }  z
\in \Omega_\delta, \ee in $\vL(-1, s; 1, -s')$.
\item More
generally, for $ 0 \le j \le n-2$ and $s > \f j 2 +1$,  there
exists some $\ep
>0$ such that \be R_0(z) = \sum_{k=0}^{j} z^{\f k 2} G_{k}  + O(|z|^{\f j 2+\ep}), \quad \mbox{ for }  z \in
\Omega_\delta, \ee in $\vL(-1, s; 1, -s)$.
\end{enumerate}
Here $G_j \in \vL(-1, s'; 1, -s')$ for $s'>\f{j}{2} +1$   and   $
0 \le j \le n-2$ are given by
\begin{eqnarray}
  G_j &= & b_0^w(v, D_x, D_v)A_j + \f{1}{k!}r_0^{(k)}(0) \quad \mbox{ for $j$ even and $j=2k$}; \\
  G_j & = & b_0^w(v, D_x, D_v)A_j \quad \mbox{ for $j$ odd}.
\end{eqnarray}
\end{lemma}

 Note that $G_j=0$ if $j$ is odd and $ 1\le j \le n-4$. The first nonzero fractional
power of $z$ in the above expansion appears for $j=n-2$. Its
coefficient $G_{n-2}$ is a rank one operator given by $ G_{n-2} =
b_0^w(v, D_x, D_v)A_{n-2} $. Using the formula
\[
\f 1 {(2\pi)^{n}}\int_{\bR_y^n} \int_{\bR^n_\xi} e^{iy\xi}
e^{2|\xi|^2} \chi(\xi) d\xi dy  = \chi(0) =1,
\]
one can compute
\be G_{n-2} =  a_{n, n-2} \w{\cdot, \vM_0}\vM_0,
\ee
 where $\vM_0$ is the free Maxwell-Boltzmann distribution:
\[
\vM_0(x, v) = (2\pi)^{-\f n 4} e^{-\f 1 4 |v|^2}, \qquad (x,v) \in \bR^{2n}.
\]
See Proposition 2.9 \cite{w2} for a similar computation when
$n=3$. The constant $a_{n, n-2}$ can be calculated by using formulas (3.4) and (3.8) in \cite{jn}. We obtain
\be \label{ann2}
a_{n, n-2} =\frac{i}{2(2 \pi)^{\frac{n-1}{2}}(n-2)!!}.
\ee

\subsection{$n\ge 4$ even.}

 When $n \ge 4$ is even, the integral kernel $(-\Delta_x -z)^{-1}$ can be
expanded asymptotically as
\begin{equation} \label{kernelneven}
A(x, y; z) = z^{\f{n-2}{2}}\ln (-i z^{\f{1}{2}}|x-y|)\sum^{+\i}_{k=0} c_{n,k} z^k|x-y|^{2k}+\f{1}{|x-y|^{n-2}} \sum^{+\i}_{k=0} d_{n,k} z^k|x-y|^{2k}.
\end{equation}
where  $c_{n,k}$ and $d_{n,k}$ are numerical constants which can be computed explicitly (\cite{jn}).  The above
expansion means that for $0 < \epsilon \le 1$, one has
\begin{align}
A(x, y; z)   =&  z^{\f{n-2}{2}}\ln (-i z^{\f{1}{2}}|x-y|)\lt( c_{n,0} + O\lt(\lt(|z||x-y|^2\rt)^\epsilon\rt)\rt)\\
             &+\f{1}{|x-y|^{n-2}} \sum^{\f{n-2}{2}}_{k=0} d_{n,k} z^k|x-y|^{2k}+ O\lt(|z|^{\f{n-2}{2}+\epsilon}|x-y|^{2\epsilon}\rt).
\end{align}
Then we have
\begin{align}
A(x, y; z)   =&  c_{n,0} z^{\f{n-2}{2}}\lt(\ln z^{\f{1}{2}}+\ln (-i|x-y|)\rt) \nn\\
              &+ z^{\f{n-2}{2}}\ln (-i z^{\f{1}{2}}|x-y|)O\lt(\lt(|z||x-y|^2\rt)^\epsilon\rt)  \label{remaider1}\\
             &+\f{1}{|x-y|^{n-2}} \sum^{\f{n-2}{2}}_{k=0} d_{n,k} z^k|x-y|^{2k}+  O\lt(|z|^{\f{n-2}{2}+\epsilon}|x-y|^{2\epsilon}\rt)\nn\\
             =& \sum^{\f{n-4}{2}}_{k=0} d_{n,k} z^k|x-y|^{2k-(n-2)}+z^{\f{n-2}{2}}\lt(c_{n,0}\ln (-i|x-y|)+ d_{n,\f{n-2}{2}}\rt) \nn\\
             &+c_{n,0} z^{\f{n-2}{2}} \ln z^{\f{1}{2}}+  O\lt(|z|^{\f{n-2}{2}+\epsilon}|x-y|^{2\epsilon}\rt) \label{remainder2}\\
              =: &  \sum^{\f{n-2}{2}}_{k=0}  z^{k} A_{2k}(x,y)+ z^{\f{n-2}{2}} \ln z^{\f{1}{2}} A_{\log}(x,y)+ A_{\text{err}}(x,y;z).\label{remainder3}
\end{align}
The remainder \eqref{remaider1} is absorbed into \eqref{remainder2}. The remainder in \eqref{remainder3} verifies:
\[
A_{\text{err}}(x,y; z)=O(|z|^{\f{n-2}{2}+\epsilon}|x-y|^{2\ep}).
\]
 $A_{\log}(x,y)=c_{n,0} $
 is a constant.  Let $A_{2k}$  denote the operator with integral kernel
 \be A_{2k}(x, y)
=  d_{n,k} |x-y|^{2k- (n-2)}, \quad k \in \bN,
\ee
for $0\le k\le \f{n-4} 2$. The continuity of $A_{2k}$ has already been studied above in the case when $n\ge 5$ is odd. The operator $A_{log}$ associated with the constant integral kernel $c_{n,0}$ is of rank one. Consequently we obtain the following

\begin{lemma} \label{lem2.6} Let $\Omega_\delta \subset \bC$ be defined as before and $n \ge 4$ even.  For  $s>\f{n}{2},\, s'>\f{n}{2}$,   there exists some $\ep>0$ such that
\be\label{evenr0}
 R_0(z) = \sum_{k=0}^{\f{n-2}{2}} z^{k} G_{2k} + z^{\f{n-2}{2}} \ln z^{\f{1}{2}} G_{\log} + O(|z|^{\f{n-2}{2}+\epsilon}), \quad \mbox{ for }  z \in
\Omega_\delta,
 \ee
 in $\vL(-1, s; 1, -s')$. Here $G_{2k} \in \vL(-1, s; 1, -s')$ for $s, s'>\f{n-4}{2}+ 2k$, $s+s'\geq 2k+2$ and $0\leq k\leq \f{n-4}{2}$, and are given by
\begin{eqnarray}
  G_{2k} &= & b_0^w(v, D_x, D_v)A_{2k} + \f{1}{k!}r_0^{(k)}(0)\\
  G_{n-2} &= & b_0^w(v, D_x, D_v)\lt(d_{n,\f{n-2}{2}}+c_{n,0}\ln(-i|x-y|\rt) \nonumber \\
   & &  + \f{1}{(\f{n-2}2)!}r_0^{(\f{n-2}2)}(0) \\
  G_{\log} &= &  b_0^w(v, D_x, D_v) A_{log}.
\end{eqnarray}
\end{lemma}

The first non integral power of $z$ in the above expansion appears for $z^{\f{n-2}2} \ln z^{\f 1 2}$ with $G_{\log}$ as its coefficient. As for the odd dimensional case, one can compute $G_{\log}$ and obtains
\be
G_{\log}= c_{n, 0} \w{\cdot, \vM_0}\vM_0.
\ee
$c_{n,0}$ can be computed by using formulas (3.4), (3.6) and (3.7) in \cite{jn}:
\be \label{cn0}
c_{n, 0}= - \frac{1}{(2 \pi)^{\frac{n}{2}}2^{\frac{n-2}{2}}\left(\frac{n-2}{2}\right)!}.
\ee

\begin{rmk} The main difference between $n\ge 5 $ odd and $n\ge 4$ even resides in the low-energy behavior of the integral kernel of $(-\Delta_x -z)^{-1}$. By (\ref{kernelneven}), in order to obtain an expansion  the form
\[
R_0(z) = G_0 + O(|z|^{\epsilon}), \quad z\in \Omega_\delta,
\]
for some $\ep >0$ in $\vL(-1, s; 1, -s')$, one already needs to take $s> \f n 2$ and $s' > \f n 2$, which are the same as in (\ref{evenr0}).
\end{rmk}

\sect{Low-energy resolvent expansions for the  full KFP operator}

\subsection{Threshold spectral properties}
Let $P = P_0 + W$ and $R(z) =(P-z)^{-1}$ for $ z\not\in \sigma(P)$.
We use the resolvent equation
\[
 R(z) =( 1 +  R_0(z)W)^{-1} R_0(z)
\]
to study the behavior of $R(z)$ near the first threshold $0$.
Assume condition (\ref{ass1}) for some $\rho >1$. Then $W \in
\vL(0, -s; -1, 1+\rho -s)$ for any $s \in \bR$. By (\ref{eq2.16})
with $k=0$, for $ \f{4-n} 2 <s'< s < \rho + \f{n-2} 2 $ with $
s-s' < \rho -1$, $G_0W \in \vL(0, -s; 1, -s')$. Since the
inclusion $\vH^{1, -s'} \hookrightarrow \vH^{-s}$ is compact when
$s> s'$, $G_0W$ is a compact operator in $\vH^{-s}$ for  $ \f{4-n}
2 < s < \rho + \f{n-2} 2 $.

\begin{theorem}\label{thm-threshold} Assume $ n\ge 4$ and $\rho >1$. Let
 $   \f{ 4-n} 2 < s < \rho + \f{ n-2} 2 $. Then
  $-1$ is not an eigenvalue of $G_0W$ in $\vH^{-s}$. In particular $ (1+ G_0W)^{-1} \in \vL(-s,
  -s)$.
\end{theorem}
\pf {\bf 1. $n \ge 5$. } Set $\delta_0 = \rho-1>0$. Let $u \in
\vH^{-s}$ with $(1+ G_0W)u =0$. We want to prove that $u $ belongs in fact to $\vH^{-s'}$ for any $s'> \f{ 4-n} 2$. Remark that $ Wu \in \vH^{-1, 2+\delta_0
-s}$. Using the equation $u = -G_0W u$ and the fact that $G_0W \in \vL(0, -s; 1, -s')$ for $s, s' \in ]\f{4-n} 2 , \rho + \f{n-2} 2 [ $ with $s-s' \le \rho -1$ ,   one obtains   $u \in \vH^{1,
-(s-\delta_0)_+}$, where for $r\in \bR$, $r_+$ is defined by
\[
r_+ =\max\{r, \f{ 4-n} 2 + \ep\}
\]
for an arbitrarily small $\ep >0$. If $s-\delta_0 \le \f{ 4-n} 2$, one has $u \in \vH^{-(\f{ 4-n} 2)_+}$.
 If $(s-\delta_0)_+ = s-\delta_0> \f{ 4-n} 2 $, we can iterate these arguments once more to
show that $u \in \vH^{2, -(s-2\delta_0)_+}$. By
iterating this argument for a finite number of times, we obtain $u
\in \vH^{2,-(\f{ 4-n} 2)_+}$.  In particular $u\in  \vH$,  if $n \ge 5$.  Note that $Pu = P_0(1+ G_0W)u=0$ which implies $u
\in D(P)$. Taking the real part of $\w{Pu,u}$, one obtains
\[
\w{(-\Delta_v + \f{|v|^2} 4 -\f n 2)u, u} =\Re \w{P u, u} =0.
\]
The positivity of operator $(-\Delta_v + \f{|v|^2} 4 -\f n 2) $ shows that $ (-\Delta_v + \f{|v|^2} 4 -\f n 2)u=0$ in $\vH$, which gives
\[
u = f(x)e^{-\f{|v|^2}4 }, \qquad \mbox{ for some }  f\in L^2_x.
\]
 Now Equation $ P u = (v\cdot \nabla_x - \nabla V(x)\cdot\nabla_v)( f(x)e^{-\f{|v|^2}4 }) =0$ gives
 $ u = c \gM$, where $\gM$ is the Maxwell-Boltzmann distribution defined by
\be \label{M-B}
 \gM(x,v) = \f{1}{ (2 \pi)^{\f n 4}} e^{- \f 12 ( \f{|v|^2}{2} + V(x))}.
 \ee
and $c$ is a constant. Since $\gM \not\in \vH$ for decaying
potential, we conclude $c=0$, hence $u=0$.

{\bf 2. $n = 4$.} For $n=4$, the argument used above only gives $
u \in \vH^{1,-s'}$ for any $s'>0$. Take $s' \in ]0, \f 12[$. Let
$\chi \in C_0^\infty(\bR)$ be a cut-off with $\chi(\tau) =1$ for $
|\tau| \le 1$ and $\chi(\tau) =0$ for $|\tau| \ge 2$ and $0 \le
\chi(\tau) \le 1$.
 Set $\chi_R(x) = \chi(\frac{|x|}{R}), R>1$ and $x\in\bR^3$  and $u_R(x,v) = \chi_R(x) u(x, v)$.
Then one has
\[
P u_R = \frac{v\cdot\hat x}{R} \chi'(\frac{|x|}{R}) u.
\]
Taking the real part of the equality
\[
\w{ P u_R, u_R}= \w{ \frac{v\cdot\hat x}{R} \chi'(\frac{|x|}{R})
u, u_R},
\]
one obtains
\be  \label{eq4.6}
\int\int_{\bR^{2n}}
|(\partial_v + \frac v 2)u(x,v)|^2 \chi(\frac{|x|}{R})^2 \; dx dv
=  \w{ \frac{v\cdot\hat x}{R} \chi'(\frac{|x|}{R}) u, u_R}
\ee
 The right hand side of  (\ref{eq4.6}) is bounded by
\[
|\w{ \frac{v\cdot\hat x}{R} \chi'(\frac{|x|}{R}) u, u_R}| \le C
R^{-(1-2s')}\|u\|_{\vH^{1,-s'}} \|u\|_{\vH^{-s'}}.
 \]
Taking the limit $R\to +\infty$ in (\ref{eq4.6}), we obtain
\be
\int\int_{\bR^{2n}} |(\partial_v + \frac v 2)u(x,v)|^2  \; dx
dv = 0.
\ee
This shows that $(\partial_v + \frac v 2)u(x,v) =0$,
a.e. in $x,v$. Since $u \in \vH^{2, -s'}$  and $Pu=0$, $u$ is of the form $u(x, v) = C(x) e^{-\frac{v^2}{4}}$ for
some $C \in L^{2,-s'}$. The arguments used above for $ n \ge
5$ allow to conclude that $u=0$.\\

Consequently  we proved for all $n\ge 4$ the injectivity of $1+
G_0W$ in $\vH^{-s}$. Fredholm Theorem for compact operators shows
that $1+ G_0W$ is invertible
 with bounded inverse in $\vH^{-s}$.
 \ef

\subsection{Low-energy resolvent expansions}

\begin{theorem}\label{thm-resol-odd} Let $n\ge 5$ be odd. Then the following results hold
true.
\begin{enumerate}
\item Assume (\ref{ass1}) with $\rho >1$ and let $s \in ]1, \f{\rho+1} 2[$. Then  there exists some constant $\delta>0$ such that for any  $z\in \Omega_\delta$, Equation $(P-z)u=0$ has no nontrivial solution in $\vH^{-s}$. In particular, $P$ has no eigenvalues in
$\Omega_\delta$. One has
 \be
R(z) = H_0 + O(|z|^\ep), \quad \mbox{ for }  z \in \Omega_\delta,
\ee
in $\vL(-1, s; 1, -s)$ for some  $\epsilon >0$. In addition, the boundary values of the resolvent
\[
R(\lambda \pm i0) =\lim_{z \to \lambda, \pm \Im z >0}  R(z),
\qquad \lambda \in [0, \delta]
\]
exist in  $\vL(-1, s; 1, -s)$ and are continuous in
$\lambda$.
\item If $\rho
> n-1$ and $s>\f n 2$, there exists some $\ep >0$
such that
\be R(z) = \sum_{k=0}^{\f{n-3} 2} z^k H_{2k} +
z^{\f{n-2} 2} H_{n-2} + O(|z|^{\f{n-2} 2+\ep}), \quad \mbox{ for }
z \in \Omega_\delta, \ee in $\vL(-1, s; 1, -s)$.
\end{enumerate}
Here $H_j\in \vL(-1, r; 1, -r)$ for $r > \f j 2 +1$. In particular, one has
\begin{eqnarray}
H_0 & = & (1 + G_0W)^{-1}G_0, \label{eqH0} \\
H_{n-2} & = &  \frac{i}{2(2 \pi)^{\frac{n-1}{2}}(n-2)!!}\w{\cdot, \gM}\gM. \label{eqHn2}
\end{eqnarray}
\end{theorem}
\pf 
We only prove Part (2), Part (1) is easier.  For
$\rho > n-1$ and $\f n 2 <s < \rho -\f{n-2} 2$, one has in
$\vL(-s; -s)$ (and also in $\vL(1,-s; 1, -s)$):
\[
1+ R_0(z) W = ( 1 + T_0(z)) (1+ G_0W)
\]
where
\[
T_0(z) = (F(z) + z^{\f{n-2} 2} G_{n-2}W )(1+ G_0W)^{-1}  + O(|z|^{\f{n-2} 2+\ep}), \mbox{ in } \vL(-s; -s)
\]
and  $ F(z) = \sum_{k=1}^{\f{n-3} 2} z^k G_{2k} W$. $T_0(z) = O(|z|)$ in $\vL(-s; -s)$. For $z\in \Omega_\delta$ with $\delta >0$ small,
$\|T_0(z)\| <1$, therefore  $1 + T_0(z) $ is invertible with uniformly bounded inverse given by the Neumann series
\[
(1+ T_0(z))^{-1} = 1 + \sum_{j=1}^\infty (-1)^j T_0(z)^j.
\]
This proves that $1 + R_0(z)W $ is invertible in $\vL(-s; -s)$ with uniformly bounded inverse given by
\begin{eqnarray*}
\lefteqn{ (1 + R_0(z)W)^{-1}}\\
 &= & (1 + G_0W)^{-1} ( 1 - T_0(z) + T_0(z)^2 + \cdots ) \\
  &=&  (1 + G_0W)^{-1} (1 + S_0(z) - z^{\f{n-2} 2} G_{n-2}W (1+ G_{n-2}W)^{-1} + O(|z|^{\f{n-2} 2+\ep} ))
\end{eqnarray*}
$z \in \Omega_\delta,$ where
\[
S_0(z) = \sum_{j=1}^{\f{n-3} 2}  z^j U_j
\]
with $U_j \in \vL(-s'; - s')$ for $s'> \f j 2 +1$. This result
implies that equation $(P-z)u =0$ has no solution in $\vH^{-s}$,
in particular, there is no eigenvalue of $P$ in $\Omega_\delta$.
From the equation $R(z) =  (1+ R_0(z) W)^{-1}R_0(z)$, it follows
that the resolvent $R(z)$  exists for $z \in\Omega_\delta$ and can
be expanded as
\begin{eqnarray*}
 R(z) &= & (1 + G_0W)^{-1} (1 + S_0(z) - z^{\f{n-2} 2} G_{n-2}W (1 + G_0W)^{-1} + O(|z|^{\f{n-2} 2+\ep} )) \\
 & & \times (\sum_{k=0}^{\f{n-3} 2} z^k G_{2k} +   z^{\f{n-2}2} G_{n-2}  + O(|z|^{\f{n-2}2 +\ep})) \\
   & = & \sum_{k=0}^{\f{n-3} 2} z^k H_{2k} +  z^{\f{n-2} 2} H_{n-2}  + O(|z|^{\f{n-2} 2+\ep}),
   \quad \mbox{ for }  z \in \Omega,\quad z \in \Omega_\delta,
\end{eqnarray*}
where $H_j \in\vL(s'; -s')$, $s'>\f j 2+1$. One can calculate: $
H_0  =  (1 + G_0W)^{-1}G_0 $ and  the coefficient  of the first
fractional power  in $z$ in this expansion is given by
\[
H_{n-2} =  (1 + G_0W)^{-1}G_{n-2} (1 - W H_0).
\]
$H_{n-2}$ is an operator of rank one in $\vL(-1, s; 1, -s)$.  As a
consequence,  the limit $R(0) = \lim_{z\in \Omega_\delta , z \to
0} R(z)$ exists and is equal to $  H_0$. Using Proposition
\ref{prop2.1}, one deduces that $R(\lambda \pm i0)$ exist in
$\vL(s', -s')$, $s'> \f 12$, for $
\lambda >0$ and small and they are continuous in $\lambda$ up to $\lambda =0$ if $s'>1$. Clearly, $R(0\pm i0) =H_0$.\\

 It remains to prove that
 \be
 (1 + G_0W)^{-1}G_{n-2} (1 - W H_0) =
a_{n,n-2} \w{\cdot, \gM} \gM,
\ee
 where $a_{n,n-2}$ is given by (\ref{ann2}). We follow the argument of
\cite{w2}.  Remark that $H_{n-2}$ is of rank one. Making use of the formula for $G_{n-2}$, one has, for $u\in
\vH^s$ with $ s > \f n 2$,
\[
H_{n-2} u =a_{n,n-2}\w{u, g_0}f_0
\]
where
\be
f_0 = (1+ G_0 W)^{-1} \gM_0, \quad  g_0 =   (1- W H_0)^*\gM_0.
 \ee
Recall that $\gM_0 \in \vH^{-s}$ for any $s> \f n 2$,  $(1 + G_0
W)^{-1} \in \vL(-s; -s)$ for any $\f n 2 <s <\rho - \f{n-2} 2$,
$G_0 W$ and $(WH_0)^{*} \in \vL( -s; -s')$ for any $s'> \f{n-4} 2$
and $\f n 2 < s  <\rho - \f{n-2} 2$. It follows that $f_0, g_0$
belong to $\vH^{-s}$ for any $s>\f n 2$.\\

 We prove first that
$f_0= \gM$.  $f_0$ satisfies the equation
\[
\lim_{z\to 0, z\not\in \bR_+} ( 1 + R_0(z) W)f_0  = (1+ G_0W) f_0 = \gM_0
\]
in $\vH^{-s}$. It follows that
\be
 Pf_0 =P_0 \gM_0 =0.
 \ee
To prove that $f_0 = \gM$, we notice that the solution of the
equation $(1+ G_0 W)\mu = \gM_0$ is unique  in $\vH^{ -s}$ for $\f n 2 <s < \rho -\f{n-2} 2$, due to the invertibility of $1+ G_0 W$. Therefore it
suffices to check that $\gM$ also verifies the equation
\be (1+
G_0 W)\gM =  \gM_0.
\ee
 To show this, remark that since $(P_0 + W)\gM =0$, one has for $\lambda <0$
\be
 ( 1+ R_0(\lambda) W \gM = - \lambda R_0(\lambda) \gM
\ee
For $\rho >\f n 2$, one has $\gM - \gM_0 \in \vH^{ s}$ for any $0\le s < \rho -\f n 2$.  By Proposition \ref{key-prop},
\[
\lim_{\lambda \to 0_-}\lambda R_0(\lambda) (\gM-\gM_0) =0, \quad \mbox{ in } \vH.
\]
Since $P_0\gM_0 =0$, $- \lambda R_0(\lambda) \gM_0 = \gM_0$ for
$\lambda <0$. It follows that
\[
(1+G_0W)\gM =  \lim_{\lambda \to 0_-}( 1+ R_0(\lambda) W) \gM =  \gM_0.
\]
The uniqueness of solution to  the equation $(1+G_0W)f_0 = \gM_0$
shows that $f_0 = \gM$. \\

To prove that $g_0=  \gM$, we remark that $1+WG_0\in \vL( s;  s)$
for any $ \f n 2 < s < \rho - \f{n-2} 2$ and it is invertible with
the inverse given by:
\be
( 1+ WG_0)^{-1} = 1 - W (1+G_0W)^{-1}G_0 = 1-W H_0.
 \ee
 Therefore
 \[
 g_0 =(1-WH_0)^*\gM_0 = (( 1+ WG_0)^{-1})^* = (1-G_0^*W)^{-1}\gM_0.
 \]
 Since $\gM$ verifies also the equation $P^*\gM =(P_0^* -W)\gM =0$, the same arguments as those used above allow to conclude that $ (1- G_0^* W)\gM = \gM_0$ which
shows
\[
g_0  = (1-G_0^*W)^{-1}\gM_0 = \gM.
\]
 This proves \eqref{eqHn2}. \ef

\begin{rmk} Making use of (\ref{A2k}), one can check that
\be
H_{j} \in \vL(-1,s'; 1, -s''), \quad \mbox{  }  0 \le j \le n-3,
\ee
for $ s' > \f{4-n} 2 +j $ and $ s''> \f{4-n} 2+ j$ with $s'+s''> 2+ j$.
\end{rmk}

\begin{theorem}\label{thm-resol-even} Let $n\ge 4$ be even. If $\rho > n-1$, $\rho+1-\f{n}{2}>s>\f n 2$, then there exists some $\ep >0$
such that \be R(z) = \sum_{k=0}^{\f{n-2} 2} z^k H_{2k} +
z^{\f{n-2} 2}\ln z^{\f{1}{2}} H_{\log} + O(|z|^{\f{n-2}{2}+\epsilon}), \quad \mbox{ for }
z \in \Omega_\delta, \ee in $\vL(-1, s; 1, -s)$.
Here $H_{2k}\in \vL(-1, r; 1, -r)$ for $r > k+ 1$. One has
\begin{eqnarray}
H_0 & = & (1 + G_0W)^{-1}G_0, \\
H_{\log} & = & - \frac{1}{(2 \pi)^{\frac{n}{2}}2^{\frac{n-2}{2}}\left(\frac{n-2}{2}\right)!} \w{\cdot, \gM}\gM. \label{eqHlog}
\end{eqnarray}
\end{theorem}

Theorem \ref{thm-resol-even} is deduced from Lemma \ref{lem2.6} and Theorem \ref{thm-threshold}, by the arguments used in the proof of Theorem \ref{thm-resol-odd}. The details are omitted here.

\sect{Large-time behavior of solutions }

The large-time behaviors of solutions to the KFP equation can be
deduced as in \cite{w2} from resolvent asymptotics and an
appropriate representation formula for the semigroup $S(t)
=e^{-tP}$. We only sketch the details here. For high energies, the
potential term $W$ can be treated as a perturbation of $P_0$. One
has

\begin{theorem}[Theorem 4.1, \cite{w2}]\label{th5.1} Let $n \ge 1$ and assume (\ref{ass1}) with $\rho \ge -1$. Then there exists
 $C>0$ such that   $\sigma (P) \cap \{z; |\Im z| > C,  \Re z \le \f 1 C |\Im z|^{\f 1 2}  \} = \emptyset$ and
\be \label{eq5.1}
\|R(z)\| \le \f{C}{|z|^{\f 1 2}}, \quad
\ee
and
\be \label{eq5.2}
\|(1-\Delta_v + v^2)^{\f 1 2}R(z)\| \le \f{C}{|z|^{\f 1 4}}, \quad
\ee
for $|\Im z| > C $ and $ \Re z \le \f 1 C |\Im z|^{\f 1 2}$.
\end{theorem}

Using this result, we can firstly represent $S(t)$ as Cauchy integral of the resolvent
\be \label{eq5.3}
S(t) f =  \frac{1}{2\pi i} \int_{\gamma } e^{-t z}R(z)f dz
\ee
for $f\in \vH$ and $t>0$, where  the contour $\gamma$ is chosen such that
\[\gamma = \gamma_- \cup \gamma_0 \cup \gamma_+\]
with $\gamma_{\pm} =\{ z; z = \pm i C + \lambda \pm i C \lambda^2,
\lambda \ge 0\}$ and $\gamma_0$ is a curve in the left-half
complex plane joining $-i C$ and $iC$ for some $C>0$ sufficiently
large, $\gamma$ being oriented from $-i\infty$ to $+i\infty$.

For $n \ge 5 $ odd and $\rho>1$, we apply  Theorem \ref{thm-resol-odd}  to deform the contour of integration to obtain
\begin{equation} \label{eq5.5}
\w{S(t)f,g} =\frac{1}{2\pi i} \int_{\Gamma} e^{-t z}\w{R(z)f,g} dz, \quad t>0,
\end{equation}
for any $f,g \in \vH^{ s}$, $s>1$. Here
\[
\Gamma = \Gamma_- \cup \Gamma_0 \cup \Gamma_+
\]
with $\Gamma_{\pm} =\{ z; z = \delta + \lambda \pm i  \delta^{-1}
\lambda^2, \lambda \ge 0\}$ for $\delta>0$ small enough and
 $\Gamma_0 =\{ z = \lambda \pm i0; \lambda \in [0, \delta]\}$. $\Gamma$ is oriented from $-i\infty$ to $+i\infty$.

\textbf{Proof of Theorem \ref{th5.4}.} For $f,g \in \vH^{s}$ with $s>1$ and for $t>0$,
 \bea \nonumber
\lefteqn{\w{S(t)f,g}} \\
 &= & \frac{1}{2\pi i} \int_{0}^{\delta} e^{- t \lambda}\w{(R(\lambda +i0)-R(\lambda -i0)) f,g} d\lambda  \\
&     + & \frac{e^{-t \delta}}{2\pi i} \int_{0}^{\infty} e^{-t (\lambda + i  \delta^{-1} \lambda^2)}\w{R(\delta + \lambda +  i  \delta^{-1} \lambda^2 ) f,g} ( 1 +  2 i \delta^{-1} \lambda)  d\lambda \nonumber  \\
 &  - &
 \frac{e^{ -t \delta}}{2\pi i} \int_{0}^{\infty} e^{-t (\lambda - i  \delta^{-1} \lambda^2)}\w{R(\delta + \lambda - i \delta^{-1} \lambda^2 ) f,g} ( 1 -  2  i \delta^{-1} \lambda)  d\lambda \nonumber
\\
& \triangleq & I_1+ I_2 + I_3. \nonumber
\eea

For $I_2$ and $I_3$, one can apply Theorem \ref{th5.1} to estimate
\[
|\w{R(\delta + \lambda \pm i  \delta^{-1} \lambda^2 ) f,g}| \le
C_M \|f\|_{\vH^{s}}\|g\|_{\vH^{s}}
\]
for $s> 1$ and $\lambda \in ]0, M]$ for each fixed $M>0$ and
\[
|\w{R(\delta + \lambda \pm i  \delta^{-1} \lambda^2 ) f,g}| \le
C_M \lambda^{-\f 1 2}\|f\|_{L^2}\|g\|_{L^2}
\]
for $\lambda >M$ with $M>1$ sufficiently large. Therefore, if
$\rho>1$ \be |I_k| \le C e^{ -t \delta_0}
\|f\|_{\vH^{s}}\|g\|_{\vH^{s}}, \quad \delta_0>0, \ee
for  $k =2,3$ and for any $s > 1$. \\

{\bf \noindent Proof of Theorem \ref{th5.4} (a).} To show
(\ref{eq5.8}), it remains  to prove that if $\rho>1$, one has
\be
\label{eq5.11} |
\int_{\Gamma_0} e^{- t z}  \w{ R(z)f,g} dz | \le C
t^{-\f n 2}\|f\|_{\vH^{s}}\|g\|_{\vH{s}}
\ee
for any $f, g \in
\vH^{s}$, $s>\f n 2$. For $ \f{4-n} 2 < s < \rho +  \f{n-2} 2$,
one has for some $\ep_0>0$
\[
W (R_0(z) - G_0) = O(|z|^{\ep_0}),  \quad \mbox{ in } \vL(s; s)
\]
for $z$ near $0$ and $z \not\in\bR_+$.  $1+ WG_0$ is invertible in
$\vL( s; s) $. It follows that
 \be
 T(z) = (1+WG_0)^{-1} W (R_0(z)
-G_0) =O(|z|^{\ep_0}), \quad \mbox{ in } \vL( s; s).
\ee
Hence $(1+ T(z))^{-1}$ can be expanded using convergent Neumann series for  $|z|$ appropriately small. Decompose $ R(z)$ as
\beas \label{eq5.12} R(z) & = & R_0(z)(1+ T(z))^{-1} (1+WG_0)^{-1} W \\
& = & \sum_{j=0}^N (-1)^j R_0(z) T(z)^j (1+WG_0)^{-1} +
O(|z|^{(N+1)\ep_0}) \eeas
 in $\vL(s; s) $, where $N$ is
taken such that $(N+1)\ep_0 > \f n 2$.  We deduce that
 \be
\label{eq5.13} |\int_{\Gamma_0} e^{- t z}  \w{ (R(z)- R_0(z)
\sum_{j=0}^N (-1)^j T(z)^j (1+WG_0)^{-1})f,g} dz |  \le C t^{-\f n 2 -
\ep}\|f\|_{\vH{s}}\|g\|_{\vH^{s}},
\ee
 for some $\ep >0$, if $ \f{4-n} 2 < s < \rho +  \f{n-2}
2$. (\ref{eq5.11}) follows from (\ref{eq5.13}) and the following lemma, which
achieves the proof of (\ref{eq5.8}).

\begin{lemma} \label{lem5.5} Let $\rho>1$. For each $j\ge 0$ and $s> \f n 2$, there exists some constant $C>0$ such that \be \label{eq5.14}
|\int_{\Gamma_0} e^{- t z}  \w{R_0(z) T(z)^j (1+WG_0)^{-1}f,g} dz
|  \le C t^{-\f n 2}\|f\|_{\vH^{s}}\|g\|_{\vH^{s}} \ee for  any
$f, g \in \vH^{s}$ and $t >1$.
\end{lemma}
\pf For $ z = \lambda + i0$ with  $\lambda \in ]0, \delta]$ and
$\delta>0$ small, \be R_0(z) - R_0(\bar z) = b_0^w(v, D_x, D_v) (
(-\Delta_x -z )^{-1} - (-\Delta_x -\bar z )^{-1}). \ee
 According to Lemma \ref{lem2.5},
  \be \label{eq5.16} R_0(z) -
R_0(\bar z)  = O( \lambda^{\f{n-2} 2}), \quad \lambda \in ]0,
\delta]
\ee
in $\vL(s;  -s)$ for any $s> \f n 2$. $ WG_0$ and
 $(1+
WG_0)^{-1}$ belong to $\vL(s'; s')$ if $ \f{4-n} 2 <s' < \rho + \f
{n-2}  2$.
 We deduce  that
 \be T(z) - T(\bar z) = (1+ WG_0)^{-1} W( R_0(z) -
R_0(\bar z)) = O( \lambda^{\f{n-2} 2}) \ee
in $\vL( s; s')$, for
$s > \f n 2$,  $ \f{4-n} 2 <s' < \rho + \f {n-2}  2$ and $s+ s'
\le 1+\rho$.  Using the formula \bea
\lefteqn{R_0(z) T(z)^j - R_0(\overline{z}) T(\overline{z})^j }\\
&= & (R_0(z)- R_0(\bar z)) T(z)^j  + R_0(\overline{z})
\sum_{k=0}^{j-1}T(z)^k (T(z)- T(\overline{z}))
T(\overline{z})^{j-k-1} \nonumber \eea one can estimate $|\w{
(R_0(z) T(z)^j - R_0(\overline{z}) T(\overline{z})^j)
(1+WG_0)^{-1}f,g}|$ by
 \bea
\lefteqn{ |\w{ (R_0(z) T(z)^j - R_0(\overline{z}) T(\overline{z})^j)  (1+WG_0)^{-1}f,g}|} \nonumber  \\[2mm]
& \le & C \{ \|(R_0(z)- R_0(\bar z))\|_{\vL( s; \; -s)} \|T(z)\|^j _{\vL(s; \; s)} \\[2mm]
  & &  +   \sum_{k=0}^{j-1} \|R_0(\overline{z})\|_{\vL( s'; \; -s)}
  \|T(z)\|^{k} _{\vL(s'; \; s')} \|(T(z)- T(\overline{z}))\|_{\vL( s; \; s')} \|T(\overline{z})\|^{j-k-1} _{\vL( s; \; s)}  \}  \nonumber  \\[2mm]
   & & \times \|f\|_{\vH^{s}}\|g\|_{\vH^{s}}  \nonumber \\[2mm]
&\le &   C' \lambda^{\f {n-2} 2} \|f\|_{\vH^{s}}\|g\|_{\vH^{s}}
\nonumber \eea for $s > \f n 2$,  $ \f{4-n} 2 <s' < \rho + \f
{n-2}  2$ and $s+ s' \le 1+\rho$. These choices are possible for
$\rho>1$.  Finally we obtain
\begin{eqnarray*}
\lefteqn{|\int_{\Gamma_0} e^{- t z}  \w{R_0(z) T(z)^j (1+WG_0)^{-1}f,g} dz |}\\
 & =&  | \int_{0}^{\delta} e^{- t \lambda} \w{ \left(R_0(\lambda+  i0) T(\lambda+  i0)^j - R_0(\lambda-  i0) T(\lambda-  i0)^j\right)  (1+WG_0)^{-1}f,g} d\lambda |\\
 &\le  & C t^{-\f n 2}\|f\|_{\vH^{s}}\|g\|_{\vH^{s}}
\end{eqnarray*}
for  any $f, g \in \vH^{s}$ and $t >1$. \ef

{\bf \noindent Proof of Theorem \ref{th5.4} (b).} Consider first $n \ge 5$ odd. Assume now $\rho>n-1$.  Then for $\f n 2 <s < \f{1+\rho} 2$, there exists some $\ep >0$ such that
\be
R(z) = \sum_{k=0}^{\f{n-3} 2} z^k H_{2k} +  z^{\f{n-2} 2} H_{n-2}  + O(|z|^{\f{n-2} 2+\ep}), \quad \mbox{ for }  z \in \Omega_\delta,
\ee
in $\vL(-1, s; 1, -s)$.
 $I_1$ can be calculated by using Theorem \ref{thm-resol-odd}
\[
I_1 = \frac{1}{\pi i} \int_{0}^{\delta} e^{- t \lambda} \w{ (
\lambda^{\f {n-2} 2}  H_{n-2} + O(\lambda^{\f 1 2 +\ep}))f,g}
d\lambda
\]
which gives \be \label{eq5.10} |I_1 - \frac{b_n}{  t^{\f n
2}}\w{H_{n-2} f,g}| \le C t^{-(\f n 2 +
\ep)}\|f\|_{\vH^{s}}\|g\|_{\vH^{s}}, \quad t \ge 1, \ee with \be
 b_n = \frac{1}{\pi i} \int_0^\infty e^{-s} s^{\f{n-2} 2} \; ds = \f 1{\pi i} \Gamma(\f n 2) = \f 1{\sqrt \pi i} \f{(n-2)!!}{2^{\f{n-1} 2}}.
\ee
Remark that
\[
a_{n,n-2} b_n = \frac{i}{2(2 \pi)^{\frac{n-1}{2}}(n-2)!!} \times  \f 1{\sqrt \pi i} \f{(n-2)!!}{2^{\f{n-1} 2}} = \f{1}{(4 \pi)^{\f n 2}}.
\]
This proves (\ref{eq5.7}) for $n \ge 5$ odd.\\

When $n \ge 4$ is even, we use Theorem \ref{thm-resol-even}:  for $\rho>n-1$ and for $\f n 2 <s < \rho+1-\f{n}{2}$,  one has
\be
R(z) = \sum_{k=0}^{\f{n-2} 2} z^k H_{2k} +  z^{\f{n-2} 2}\ln z^{\f{1}{2}} H_{\log}  + O(|z|^{\f{n-2} 2+\epsilon}), \quad \mbox{ for }  z \in \Omega_\delta,
\ee
in $\vL(-1, s; 1, -s)$.
 $I_1$ can be calculated by
\begin{align}
I_1 =& \frac{1}{2\pi i} \int_{0}^{\delta} e^{- t \lambda} \lambda^{\f {n-2} 2}\w{ (
\lt(\ln (\la+i0)^{1/2}-\ln(\la-i0)^{1/2}\rt)  H_{\log} + O(\lambda^{\f{ n-2}{ 2} +\ep}))f,g}
d\lambda\nn\\
=& \frac{1}{2\pi i} \int_{0}^{\delta} e^{- t \lambda} \lambda^{\f {n-2} 2}\w{ -i\pi H_{\log} f,g}
d\lambda+\frac{1}{2\pi i} \int_{0}^{\delta} e^{- t \lambda} \w{ O(\lambda^{\f{ n-2}{ 2} +\ep}))f,g}
d\lambda
\end{align}
which gives \be \label{eq5.10b} |I_1 - \frac{e_n}{  t^{\f n
2}}\w{H_{\log} f,g}| \le C t^{-(\f n 2 +
\ep)}\|f\|_{\vH^{s}}\|g\|_{\vH^{s}}, \quad t \ge 1, \ee with \be
 e_n = -\frac{1}{2} \int_0^\infty e^{-s} s^{\f{n-2} 2} \; ds = - \f 1 2 \left( \f{n-2} 2 \right)!.
\ee
Noticing that
\[
c_{n,0} e_n = \f{1}{ ( 4 \pi)^{\f n 2}},
\]
 (\ref{eq5.7})  for $ n\ge 4$ even follows.
\ef

\sect{Dispersive estimates for quickly decaying potentials}

One can derive from Theorem \ref{th5.4} global-in-time $L^p-L^q$ estimates for the semi-group $ S(t) =e^{-tP}$, $t>0$, where
\[
L^p = L^p(\bR^{2n}_{x,v}; dx dv)
\]
is equipped with the natural norm,   $1 \le p < q \le \i$ and $\f 1 p + \f 1 q =1$. Under fairly general conditions, $e^{-tP}$ is a strongly continuous positivity preserving contraction semigroup in $L^p$. It can be extended as a map from $L^p$ to $L^q$. See \cite{wz}.

\begin{theorem}\label{th6.1} Assume condition (\ref{ass1}) with $\rho >1$ if $n \ge 5$ is odd, and $\rho >n-1$ if $n \ge 4 $ is even.  Then for $1 \le p < q \le \i$ and $\f 1 p + \f 1 q =1$, there exists some constant $C>0$ such that
\be\label{eq6.1}
\|e^{-tP}\|_{\vL(L^p; L^q)} \le \f{C}{(\gamma(t))^{\f n {2 p}(1-\f p q)}}, \quad t\in ]0, \i[,
\ee
 where $\gamma(t) = \sigma(t) \theta(t) $ with
\be \label{eq6.2}
  \quad \sigma(t) = t - 2 \coth (t) + 2 \text{\rm cosech}(t),  \quad  \theta(t) = 4\pi  e^{- t } \sinh (t).
\ee
\end{theorem}

The function $\gamma(t)$ appears in the fundamental solution to the free KFP equation (\cite{wz}). It behaves like
\be
\gamma(t) = \left\{ \begin{array}{ccc} \f \pi 3 t^4, & & \mbox { as } t \to 0 \\
 4\pi t, & & \mbox{ as } t \to +\infty\end{array}
 \right.
\ee
The short-time estimates in \eqref{eq6.1} are conform with the global subellipticity  of the KFP operator with loss of $\f 1 3$ derivatives, and the large-time ones can be compared with the dispersive estimates  for the heat equation. \\

The proof of Theorem \ref{th6.1} follows the same argument of \cite{wz}. In fact, the estimate \eqref{eq6.1} for $t\in ]0, 1]$ is already proven in Theorem 4.1 \cite{wz} which is valid for $\rho \ge -1$ and $n\ge 1$. For  $n\ge 5$ and $\rho >1$, we can use Duhamel's principle and Theorem \ref{th5.4} (a) to prove \eqref{eq6.1} for $t\ge 1$.  When $n\ge 4$ is even, we need  the condition \eqref{ass1} with  $\rho > n-1$ in order to apply Theorem \ref{th5.4} (b). The details are the same as in \cite{wz} for $n=3$ and are omitted here.

\section*{Acknowledgments}

{ X. Pan is supported by NSF of China  under Grant  No.12471222 and No.12031006. }

\end{document}